\documentclass{amsart}
\usepackage{graphicx}
\usepackage[top=1.5in,bottom=1.5in,left=1.5in,right=1.5in]{geometry}
\usepackage{leftidx}
\newtheorem{theorem}{Theorem}[section]
\newtheorem*{theorem A}{Theorem A}
\newtheorem*{theorem B}{N\"olker's Theorem}
\newtheorem{lemma}{Lemma}[section]

\newtheorem{proposition}{Proposition}[section]

\newtheorem{corollary}{Corollary}[section]

\theoremstyle{remark}

\theoremstyle{remark}

\theoremstyle{definition}

\numberwithin{equation}{section}
\def\({\left ( }
\def\){\right )}
\def\<{\left < }
\def\>{\right >}

\begin{document}

\title{Real hypersurfaces with $^{*}$-Ricci solitons of non-flat complex space forms}

\author{Xiaomin Chen}
\address{College of  Science, China University of Petroleum-Beijing, Beijing, 102249, China}
\email{xmchen@cup.edu.cn}

\thanks{The author is supported by the Science Foundation of China
University of Petroleum-Beijing(No.2462015YQ0604) and partially
by  the Personnel Training and Academic
Development Fund (2462015QZDX02).}

%SRFDP of China (Grant No. 20050027025).

\begin{abstract}
Kaimakamis and Panagiotidou in \cite{KP} introduced the notion of $^*$-Ricci soliton
and studied the real hypersurfaces of a non-flat complex space form admitting a $^*$-Ricci soliton whose potential vector
field is the structure vector field.  In this article, we consider that a real hypersurface of a non-flat complex
space form admits a $^*$-Ricci soliton whose potential vector field belongs to the principal curvature space
and the holomorphic distribution.
\end{abstract}



\maketitle

\section{Introduction}

An $n$-dimensional complex space form is an $n$-dimensional K\"ahler manifold with constant sectional
curvature $c$. A complete and simple connected complex space form with $c\neq0$(i.e. a complex projective space
 $\mathbb{C}P^n$ or a complex hyperbolic space $\mathbb{C}H^n$)
 is called a \emph{non-flat complex space form} and denoted by $\widetilde{M}^n(c)$.

 Let $M$ be a real hypersurface of $\widetilde{M}^n(c)$, then
 there exists an almost contact structure $(\phi,\eta,\xi,g)$ on $M$ induced
 from $\widetilde{M}^n(c)$. The study of real hypersurfaces
in a non-flat complex space form is a very interesting and active field in recent decades
and many results of the classification of real hypersurfaces in non-flat complex space
forms were achieved (cf.\cite{B,K2,T1,T2,W}). In particular, if $\xi$ is an eigenvector of shape operator $A$ then $M$ is
called a \emph{Hopf hypersurface}, and
we note that the following conclusion is due to Kimura and Takagi for $\mathbb{C}P^n$ and Berndt for $\mathbb{C}H^n$.
\begin{theorem}[\cite{KI},\cite{T3},\cite{B}]\label{A}
Let $M$ be a Hopf hypersurface in non-flat complex space form $\widetilde{M}^n(c),n\geq2$. If $M$ has constant principal curvatures,
then the classification is as follows:\\
$\bullet$ In case of $\mathbb{C}P^n$, $M$ is locally congruent to one of the following:
\begin{enumerate}
  \item $A_1$: Geodesic hyperspheres.
  \item $A_2$: Tubes over a totally geodesic complex projective space $\mathbb{C}P^k$ for $1\leq k\leq n-2$.
\item $B$:  Tubes over a complex quadric $Q_{n-1}$ and $\mathbb{R}P^n$.
\item $C$: Tubes over Segre embedding of $\mathbb{C}P^1\times\mathbb{C}P^{\frac{n-1}{2}}$ and $n(\geq5)$ and is odd.
\item $D$: Tubes over Pl\"ucker embedding of the complex Grassmannian manifold $G_{2,5}$. Occur only for $n=9$.
\item $E$:  Tubes over the canonical embedding Hermitian symmetry space $SO(10)/U(5)$. Occur only for $n=15$.
\end{enumerate}

$\bullet$ In case of $\mathbb{C}H^n$, $M$ is locally congruent to one of the following:
\begin{enumerate}
 \item $A_1$: Geodesic hyperspheres $(Type\, A_{11})$ and tubes over totally geodesic complex hyperbolic hyperplanes $(Type\, A_{12})$.
 \item $A_2$:  Tubes over totally geodesic $\mathbb{C}H^{k}\subset\mathbb{C}H^n$ for some $k\in\{1,\cdots,n-2\}$.
  \item $B$: Tubes over a totally geodesic real hyperbolic space $\mathbb{R}H^n\subset\mathbb{C}H^n$.
\item $N$: Horospheres.
\end{enumerate}
\end{theorem}

In particular, if $M$ has two distinct constant principal curvatures, the classification is as follows:
\begin{theorem}[\cite{T1}, Corollary 2 in \cite{BD2}]\label{B}
Let $M$ be a hypersurface in non-flat complex space form $\widetilde{M}^n(c),n\geq2$, with two distinct constant principal curvatures. Then\\
$\bullet$ in case of $\mathbb{C}P^n$,
$M$ is locally congruent geodesic hyperspheres in $\mathbb{C}P^n$(Type $A_1$);\\
$\bullet$ in case of $\mathbb{C}H^n$, $M$ is locally congruent to one of the following:
\begin{enumerate}
\item $A_{11}$: Geodesic hyperspheres in $\mathbb{C}H^n$.
 \item $A_2$: Tubes around a totally geodesic $\mathbb{C}H^{n-1}\subset\mathbb{C}H^n$.
  \item $B$: Tubes of radius $r=\ln(2+\sqrt{3})$ around a totally geodesic real hyperbolic space $\mathbb{R}H^n\subset\mathbb{C}H^n$.
\item $N$: Horospheres in $\mathbb{C}H^n$.
 \end{enumerate}
\end{theorem}

Since there are no Einstein real hypersurfaces in $\widetilde{M}^n(c)$ (see \cite{CR} and \cite{M}), Cho and Kimura in
\cite{CK} considered a real hypersurface in $\widetilde{M}^n(c)$ admitting a Ricci soliton. The notion of Ricci soliton,
introduced firstly by Hamilton in \cite{H}, is the generalization
of Einstein metric, that is, a Riemannian metric $g$ satisfying
\begin{equation*}
 \frac{1}{2}\mathcal{L}_W g+Ric-\lambda g=0,
\end{equation*}
where $\lambda$ is a constant and $\mathrm{Ric}$ is the Ricci tensor of $M$. The vector field $W$ is
called \emph{potential vector field}.  Moreover, the Ricci soliton is called
 shrinking, steady and expanding according as $\lambda$ is positive, zero and negative respectively.
 In \cite{CK}, it is proved that there does not admit a Ricci soliton on $M$ when the potential vector field is the structure field $\xi$.
At the same time, by introducing a so-called $\eta$-Ricci soliton $(\eta,g)$ on $M$, which satisfies
\begin{equation*}
 \frac{1}{2}\mathcal{L}_W g+Ric-\lambda g-\mu \eta\otimes\eta=0,
\end{equation*}
for constants $\lambda,\mu$, they gave a classification of a real hypersurface admitting an $\eta$-Ricci soliton
whose potential vector is the structure field $\xi$. In \cite{CK2}, Cho and Kimura also proved that
 a compact real hypersurface of contact-type in a complex number space admitting a Ricci soliton is a sphere and
 a compact Hopf hypersurface in a non-flat complex space form does not admit a Ricci soliton.

As the corresponding of Ricci tensor, in \cite{HA} Hamada  defined the $^*$-Ricci tensor $\mathrm{Ric^*}$ in real hypersurfaces of complex space form as
\begin{equation*}
  Ric^*(X,Y)=\frac{1}{2}(trace\{\phi\circ R(X,\phi Y)\}),\quad\hbox{for all}\; X,Y\in TM,
\end{equation*}
and if the $^*$-Ricci tensor is a constant multiple of $g(X,Y)$ for all $X,Y$ orthogonal to $\xi$, then
$M$ is said to be a \emph{$^*$-Einstein manifold}.  Furthermore, Hamada gave the following result of the $^*$-Einstein Hopf hypersurfaces in non-flat space forms.
\begin{theorem}[\cite{HA}]\label{C}
Let $M$ be a $^*$-Einstein Hopf hypersurface in non-flat complex space form $\widetilde{M}^n(c),n\geq2$.\\
$\bullet$ In case of $\mathbb{C}P^n$, $M$ is an open part of one of the following:
\begin{enumerate}
  \item $A_1$: a geodesic hypersphere;
  \item $A_2$: a tuber over a totally geodesic complex projective space $\mathbb{C}P^k$ of radius $\frac{\pi r}{4}$ for $1\leq k\leq n-2$, where $r=\frac{2}{\sqrt{c}}$;
  \item $B$: a tuber over a complex quadric $Q_{n-1}$ and $\mathbb{R}P^n$.
\end{enumerate}
$\bullet$ In case of $\mathbb{C}H^n$, $M$ is an open part of one of the following
\begin{enumerate}
  \item $A_{11}$: a geodesic hypersphere;
  \item $A_{12}:$ a tube around a totally geodesic complex hyperbolic hyperplane;
  \item $B$: a tube around a totally geodesic real hyperbolic space $\mathbb{R}H^n$;
  \item $N$: a horosphere.
 \end{enumerate}
\end{theorem}

Motivated by the work in \cite{CK,CK2,HA}, Kaimakamis and Panagiotidou in \cite{KP} introduced a so-called $^*$-Ricci soliton, that is, a
Riemannain metric $g$ on $M$ satisfying
\begin{equation}\label{1}
 \frac{1}{2}\mathcal{L}_W g+Ric^*-\lambda g=0,
\end{equation}
where $\lambda$ is constant and $\mathrm{Ric^*}$ is the $^*$-Ricci tensor of $M$.
 They considered the case where $W$ is the structure field $\xi$ and obtained that a real hypersurface in complex projective space does
not admit a $^*$-Ricci soliton and a real hypersurface in complex hyperbolic space admitting a $^*$-Ricci soltion is locally congruent to a geodesic hypersphere.

It is well-known that the tangent bundle $TM$ can be decomposed as $TM=\mathbb{R}\xi\oplus\mathcal{D}$, where $\mathcal{D}=\{X\in TM,\eta(X)=0\}$ is called \emph{holomorphic distribution}. In the last part of \cite{KP}, they proposed two open problems:
\begin{description}
  \item[Problem 1] Are there real hypersurfaces admitting a $^*$-Ricci soliton whose potential vector field is a principal vector field of the real hypersurface?
  \item[Problem 2] Are there real hypersurfaces admitting a $^*$-Ricci soliton whose potential vector field belongs to the holomorphic distribution $\mathcal{D}$?
\end{description}

In the present paper, we shall consider the above two problems. For the Problem 1, we consider the case of 2-dimensional non-flat complex space forms.
Denote by $T_\chi$ the distribution on $M$ formed by principal curvature spaces of $\chi$ and $\Gamma(T_\chi)$
by the all smooth sections of $T_\chi$. We obtain the following conclusions:
\begin{theorem}\label{T1}
Let $M$ be a hypersurface of non-flat complex space form $\widetilde{M}^2(c)$ with a $^*$-Ricci soliton whose
potential vector field $W\in\Gamma(T_\chi),\chi\neq0$. If the principal curvatures are  constant along $\xi$ and $A\xi$
 then\\
$\bullet$ in case of $\mathbb{C}P^2$, $M$ is an open part of a tube around the complex quadric, or a geodesic hypersphere;\\
$\bullet$ in case of $\mathbb{C}H^2$, $M$ is an open part of
\begin{enumerate}
  \item [(1)] a geodesic hypersphere, or
  \item[(2)] a tubes around a totally geodesic $\mathbb{C}H^1$, or
  \item[(3)] a tubes around a totally geodesic real hyperbolic space $\mathbb{R}H^2$, or
   \item[(4)] a horosphere.
 \end{enumerate}
\end{theorem}
\begin{theorem}\label{T2}
Let $M$ be a hypersurface of complex projective space $\mathbb{C}P^2$, admitting a $^*$-Ricci soliton whose
potential vector field $W\in\Gamma(T_0)$. Then $M$ is an open part of a tube around the complex quadric.
\end{theorem}

For the Problem 2, we first obtain the following result:
\begin{theorem}\label{T3}
Let $M$ be a hypersurface of complex projective space $\mathbb{C}P^2$ with a $^*$-Ricci soliton whose
potential vector field $W\in\mathcal{D}$.  If the principal curvatures are constant along $\xi$ and $A\xi$,
 then $M$ is locally congruent to a geodesic hypersphere in $\mathbb{C}P^2$. Moreover, if $g(A\xi,\xi)=0$ then $W$ is Killing.
\end{theorem}

Furthermore, due to the decomposition $TM=\mathbb{R}\xi\oplus\mathcal{D}$, we have $A\xi=a\xi+V,$
where $V\in\mathcal{D}$ and $a$ is a smooth function on $M$. The following conclusion is obtained:
\begin{theorem}\label{T4}
Let $M^{2n-1}$ be a hypersurface of complex space form $\widetilde{M}^n(c)$, $n\geq2$. Then \\
$\bullet$ in case of $\mathbb{C}P^n$ there are no real hypersurfaces admitting a $^*$-Ricci soliton with
potential vector field $W=V$;\\
$\bullet$ in case of $\mathbb{C}H^n$, if $M$ admits a $^*$-Ricci soliton with potential vector field $W=V$, it is locally congruent to a geodesic hypersphere.
\end{theorem}
This paper is organized as follows. In Section 2, some basic concepts and formulas
are presented. To prove $M$ is Hopf under the assumptions of theorems, in Section 3 we give some formulas for the non-Hopf hypersurfaces with $^*$-Ricci solitons, and the proofs of theorems are given in Section 4, Section 5 and Section 6, respectively.

\section{Preliminaries}

  Let ($\widetilde{M}^n,\widetilde{g})$ be a complex $n$-dimensional K\"ahler manifold
   and $M$ be an immersed real hypersurface of $\widetilde{M}^n$ with induced metric $g$.
We denote by $J$ the complex structure on $\widetilde{M}^n$. There exists a local defined
unit normal vector field $N$ on $M$ and we write $\xi:=- JN$
by the structure vector field of $M$.
 An induced one-form $\eta$ is defined by
$\eta(\cdot)=\widetilde{g}(J\cdot,N)$, which is dual to $\xi$.  For any vector field $X$ on $M$ the tangent part of $JX$
is denoted by $\phi X=JX-\eta(X)N$. Moreover, the following identities hold:
\begin{equation}\label{2}
\phi^2=-Id+\eta\otimes\xi,\quad\eta\circ \phi=0,\quad\phi\circ\xi=0,\quad\eta(\xi)=1,
\end{equation}
\begin{equation}\label{3}
g(\phi X,\phi Y)=g(X,Y)-\eta(X)\eta(Y),
\end{equation}
\begin{equation}\label{4}
g(X,\xi)=\eta(X),
\end{equation}
where $X,Y\in\mathfrak{X}(M)$. By \eqref{2}-\eqref{4}, we know that $(\phi,\eta,\xi,g)$ is an almost
contact metric structure on $M$.

Denote by $\nabla, A$ the induced Riemannian connection and the shape operator on $M$, respectively.
Then the Gauss and Weingarten formulas are given by
\begin{equation}\label{5}
\widetilde{\nabla}_XY=\nabla_XY+g(AX,Y)N,\quad\widetilde{\nabla}_XN=-AX,
\end{equation}
 where $\widetilde{\nabla}$ is the connection on $\widetilde{M}^n$ with respect to $\widetilde{g}$.
Also, we have
\begin{equation}\label{6}
  (\nabla_X\phi)Y=\eta(Y)AX-g(AX,Y)\xi,\quad\nabla_X\xi=\phi AX.
\end{equation}
$M$ is said to be a \emph{Hopf hypersurface} if the structure vector field $\xi$ is an eigenvector of $A$.

From now on we always assume that the sectional curvature of $\widetilde{M}^n$ is contant $c\neq0$, i.e. $\widetilde{M}^n$ is a non-flat complex space form, denoted by $\widetilde{M}^n(c)$, then the curvature tensor $R$ of $M$ is given by
 \begin{equation}\label{7}
 \begin{aligned}
R(X,Y)Z=\frac{c}{4}\Big(&g(Y,Z)X-g(X,Z)Y+g(\phi Y,Z)\phi X
-g(\phi X,Z)\phi Y\\
&+2g(X,\phi Y)\phi Z)\Big)+g(AY,Z)AX-g(AX,Z)AY,
\end{aligned}
\end{equation}
and the shape operator $A$ satisfies
\begin{equation}\label{8}
 (\nabla_XA)Y-(\nabla_YA)X=\frac{c}{4}\Big(\eta(X)\phi Y-\eta(Y)\phi X-2g(\phi X,Y)\xi\Big),
\end{equation}
for any vector fields $X,Y,Z$ on $M$.

Recall that the $^*$-Ricci operator $Q^*$ of $M$ is defined by
\begin{equation*}
  g(Q^*X,Y)=Ric^*(X,Y)=\frac{1}{2}trace\{\phi\circ R(X,\phi Y)\},\quad\hbox{for all}\, X,Y\in TM.
\end{equation*}
By \eqref{7}, it is proved in Theorem 2 of \cite{IR} that the $^*$-Ricci operator is expressed as
\begin{equation}\label{9}
Q^*=-\Big[\frac{cn}{2}\phi^2+(\phi A)^2\Big].
\end{equation}
In particular, if $Q^*=0$ then $M$ is said to be a \emph{$^*$-Ricci flat hypersurface}.
Due to \eqref{2} the $^*$-Ricci soliton equation \eqref{1} becomes
\begin{equation}\label{10}
\begin{aligned}
  g(\nabla_XW,Y)+g(X,\nabla_YW)&+ncg(X,Y)-nc\eta(X)\eta(Y)\\
   &+2g(\phi AX,A\phi Y)-2\lambda g(X,Y)=0,
\end{aligned}
\end{equation}
for any vector fields $X,Y$ on $M$.

\section{Non-Hopf hypersurfaces with $^*$-Ricci solitons}
In this section we assume that $M$ is a non-Hopf hypersurface in $\widetilde{M}^2(c)$ with a $^*$-Ricci soliton. Since $M$ is not Hopf,
 due to the decomposition $TM=\mathbb{R}\xi\oplus\mathcal{D}$,
we can write $A\xi$ as
\begin{equation}\label{11}
  A\xi=\alpha\xi+\beta U,
\end{equation}
where $\alpha=\eta(A\xi),\beta=|\phi\nabla_\xi\xi|$ are the smooth functions on $M$ and
 $U=-\frac{1}{\beta}\phi\nabla_\xi\xi\in\mathcal{D}$ is a unit vector field with $\beta\neq0$.  Write
 $$\mathcal{N}:=\{p\in M:\beta\neq0\quad\hbox{in a neighbourhood of}\;p\}.$$
\begin{lemma}\label{L4.1}
On $\mathcal{N}$, we have $A\phi U=0.$
\end{lemma}
\proof In view of $^*$-Ricci soliton equation \eqref{1}, we know $Ric^*(X,Y)=Ric^*(Y,X)$ for every vector fields $X,Y\in TM$.
That means that for every vector field $X$,
\begin{equation}\label{12}
  \phi A\phi AX=A\phi A\phi X.
\end{equation}
On the other hand, we have
\begin{align*}
  \phi^2 A\phi AX &=-A\phi AX+\eta(A\phi AX)\xi \\
   & =-A\phi AX+g(\alpha\xi+\beta U,\phi AX)\xi\\
&=-A\phi AX-\beta g(\phi U,AX)\xi
\end{align*}
and
\begin{align*}
  \phi A\phi A\phi X &=A\phi A\phi^2 X \\
   & =-A\phi AX+\eta(X)A\phi A\xi\\
&=-A\phi AX+\beta\eta(X)A\phi U.
\end{align*}
Since $\beta\neq0$ on $\mathcal{N}$, we get from \eqref{12} that $-g(\phi U,AX)\xi=\eta(X)A\phi U$. Taking $X=\xi$ in this formula we obtain the desired result.\qed\bigskip

Since $\{\xi,U,\phi U\}$ is a locally orthonormal frame on $\mathcal{N}$, there are smooth functions $\gamma,\mu,\delta$ such that
\begin{equation}\label{13}
  AU=\beta\xi+\gamma U+\delta\phi U,\qquad A\phi U=\delta U+\mu\phi U.
\end{equation}
By Lemma \ref{L4.1}, we have $\delta=\mu=0.$
Moreover, in \cite{PX} the following lemma was proved:
\begin{lemma}\label{L2}
With respect to the orthonormal basis $\{\xi,U,\phi U\}$, we have
\begin{equation*}\begin{aligned}
  &\nabla_U\xi=\gamma\phi U,\quad\quad\nabla_{\phi U}\xi=0,\quad\nabla_\xi\xi=\beta\phi U, \\
 &\nabla_UU=k_1\phi U,\quad\quad\nabla_{\phi U}U=k_2\phi U,\quad\nabla_\xi U=k_3\phi U, \\
 &\nabla_U\phi U=-k_1U-\gamma\xi,\quad\nabla_{\phi U}\phi U=-k_2U,\quad\nabla_\xi\phi U=-k_3U-\beta\xi,
\end{aligned}
\end{equation*}
where $k_1,k_2,k_3$ are smooth functions on $M$.
\end{lemma}

Applying Lemma \ref{L2}, we have
\begin{proposition} The following formulas on $\mathcal{N}$ are valid:
\begin{align}
 & k_3\beta+\alpha\beta-\phi U(\alpha) =0,\quad k_2=0,\label{14} \\
  &k_3\gamma+\beta^2-\phi U(\beta) =-\frac{c}{4}\label{15},\\
  &\xi(\beta)=U(\alpha),\quad\xi(\gamma)=U(\beta),\label{16} \\
  &\beta^2+k_3\gamma-\alpha\gamma-\beta k_1=\frac{c}{4}\label{17},\\
  &k_1\beta+\alpha\gamma-\phi U(\beta) =-\frac{c}{2}. \label{18}
\end{align}
\end{proposition}
\proof  By taking $X=\xi$ and $Y=\phi U$ in \eqref{8}, we obtain
\begin{equation*}
  (\nabla_\xi A)\phi U-(\nabla_{\phi U}A)\xi=-\frac{c}{4}U.
\end{equation*}
In view of \eqref{13} and Lemma \ref{L2}, the above formula leads to $k_2=0$ since $\beta\neq0$.
Also \eqref{14} and \eqref{15} are attained.
By a straightforward computation, the relation \eqref{8} for $X=\xi$ and $Y=U$ implies \eqref{16} and \eqref{17}.
And the relation \eqref{8} for $X=U$ and $Y=\phi U$ gives \eqref{18}.\qed\bigskip

Let us assume that $W$ is an eigenvector of $A$, namely, there is a smooth function $\chi$ such that $AW=\chi W$ holds. On $\mathcal{N}$,
in the basis of $\{\xi,U,\phi U\}$ the potential vector $W$ may be expressed as
\begin{equation*}
  W=f_1\xi+f_2U+f_3\phi U,
\end{equation*}
where $f_1,f_2,f_3$ are the smooth functions on $\mathcal{N}$.

In view of Lemma \ref{L2}, by a direct computation, we have
\begin{align}
\nabla_\xi W&=(\xi(f_1)-f_3\beta)\xi+(\xi(f_2)-f_3k_3)U+(f_1\beta+f_2k_3+\xi(f_3))\phi U,\label{19}\\
\nabla_UW&=(U(f_1)-f_3\gamma)\xi+(U(f_2)-f_3k_1)U+(f_1\gamma+f_2k_1+U(f_3))\phi U,\label{20}\\
\nabla_{\phi U}W&=\phi U(f_1)\xi+\phi U(f_2)U+\phi U(f_3)\phi U.\label{21}
\end{align}

Putting $X=Y=\xi$ in \eqref{10}, by \eqref{19} we find
\begin{equation}\label{22}
  \xi(f_1)-f_3\beta=\lambda.
\end{equation}
Furthermore, putting $X=Y=U$ and $X=Y=\phi U$ in \eqref{10} respectively,  we get from \eqref{20} and \eqref{21} that
\begin{align}
  U(f_2)-f_3k_3+c-\lambda=0,\label{23} \\
  \phi U(f_3)+c-\lambda=0.\label{24}
\end{align}
Also, when $X$ and $Y$ are taken as the different vectors of $\xi,U$, and $\phi U$ in \eqref{10}, a similar computation leads to
\begin{equation}\label{25}
\left\{
  \begin{array}{ll}
  \xi(f_2)-f_3k_3+U(f_1)-f_3\gamma =0, \\
  f_1\beta+f_2k_3+\xi(f_3)+\phi U(f_1)=0,\\
f_1\gamma+f_2k_1+U(f_3)+\phi U(f_2)=0.
  \end{array}
\right.
\end{equation}

Actually, Lemma \ref{L4.1} shows that at every point of $\mathcal{N}$ there exists a principal curvature $0$ and $\phi U$ is the corresponding principal vector.
It turns out that there are at least two distinct principal curvatures in non-flat complex space forms (see \cite[Theorem 1.5]{NR}).

Let $\lambda_i$ be the principal curvatures for $i=1,2,3$, where $\lambda_3=0$.
We may assume that $e_1=\cos\theta\xi+\sin\theta U,e_2=\sin\theta\xi-\cos\theta U$ are the unit principal vectors corresponding to $\lambda_1$ and $\lambda_2$, respectively, where $\theta$ is the angle between principal vector $e_1$ and $\xi$.
It is clear that $\{e_1,e_2,e_3=\phi U\}$ is also an orthonormal frame.
Namely,
\begin{equation*}
A(e_1,e_2,e_3)=(e_1,e_2,e_3)\left(
  \begin{array}{ccc}
    \lambda_1 &  &  \\
     & \lambda_2 &  \\
     &  & 0 \\
  \end{array}
\right).
\end{equation*}
Denote by
\begin{equation*}
B=\left(
    \begin{array}{ccc}
      \cos\theta & \sin\theta & 0 \\
      \sin\theta & -\cos\theta & 0 \\
      0 & 0 & 1 \\
    \end{array}
  \right)
\end{equation*}
the transformation matrix of two frames, i.e.
\begin{equation*}
  (e_1,e_2,e_3)=(\xi,U,\phi U)B.
\end{equation*}
Moreover, since
\begin{equation*}
  A(\xi,U,\phi U)=(\xi,U,\phi U)\left(
     \begin{array}{ccc}
       \alpha & \beta & 0 \\
       \beta & \gamma & 0 \\
       0 & 0 & 0 \\
     \end{array}
   \right),
\end{equation*}
we get
\begin{equation*}
 \left(
     \begin{array}{ccc}
       \alpha & \beta & 0 \\
       \beta & \gamma & 0 \\
       0 & 0 & 0 \\
     \end{array}
   \right)=
B\left(
  \begin{array}{ccc}
    \lambda_1 &  &  \\
     & \lambda_2 &  \\
     &  & 0 \\
  \end{array}
\right)B^T.
\end{equation*}
A straightforward calculation leads to
\begin{equation}\label{26}
    \alpha=\lambda_1\cos^2\theta+\lambda_2\sin^2\theta,\quad\beta=\frac{1}{2}(\lambda_1-\lambda_2)\sin2\theta,
\quad\gamma=\lambda_1\sin^2\theta+\lambda_2\cos^2\theta.
\end{equation}

If $M$ has only two distinct  principal curvatures at any point $p\in\mathcal{N}$, then either $\lambda_1=\lambda_2\neq0$, or one of $\lambda_1$ and $\lambda_2$ vanishes. However, the second of \eqref{26} will come to $\beta=0$ if $\lambda_1=\lambda_2$, thus it is impossible.
Without loss generality, we set $\lambda_1=0$ and $\lambda_2\neq0$.
In terms of \cite[Theorem 4]{IR2}, $\alpha,\beta$ and $\gamma$ satisfy
\begin{align*}
\xi(\alpha)=\xi(\beta)&=\xi(\gamma)=0,\\
  U(\alpha)=&\beta(\alpha+\gamma).
\end{align*}
Using \eqref{16}, we thus derive $\alpha+\gamma=0$ because $\beta\neq0$. This shows $\lambda_2=0$ from the first and third of \eqref{26}. It is a contradiction.
Therefore on $\mathcal{N}$ there are three distinct principal curvatures, i.e. $\lambda_1,\lambda_2$ are not zero and $\lambda_1\neq\lambda_2$.

Using \eqref{16} again, we derive from \eqref{26} that
\begin{align*}
  U(\lambda_1)\cos^2\theta+U(\lambda_2)\sin^2\theta-(\lambda_1-\lambda_2)\sin2\theta U(\theta)\\
=\frac{1}{2}\xi(\lambda_1-\lambda_2)\sin2\theta+(\lambda_1-\lambda_2)\cos2\theta \xi(\theta),\\
 \xi(\lambda_1)\sin^2\theta+\xi(\lambda_2)\cos^2\theta+(\lambda_1-\lambda_2)\sin2\theta \xi(\theta)\\
=\frac{1}{2}U(\lambda_1-\lambda_2)\sin2\theta+(\lambda_1-\lambda_2)\cos2\theta U(\theta).
\end{align*}
From which we arrive at
\begin{align*}
  \xi(\theta)=\frac{U(\lambda_1-\lambda_2)+U(\lambda_1+\lambda_2)\cos2\theta-\xi(\lambda_1+\lambda_2)\sin2\theta}{2(\lambda_1-\lambda_2)},\\
U(\theta)=\frac{-\xi(\lambda_1-\lambda_2)+\xi(\lambda_1+\lambda_2)\cos2\theta+U(\lambda_1+\lambda_2)\sin2\theta}{2(\lambda_1-\lambda_2)}.
\end{align*}
Thus we obtain
\begin{proposition}\label{P3.2}
 If on $\mathcal{N}$ the principal curvatures are constant along $\xi$ and $A\xi$, then the following equations hold:
\begin{align}
  &\xi(\theta)=U(\theta)=0,\label{27}\\
&\xi(\beta)=U(\alpha)=\xi(\gamma)=U(\beta)=0.\label{28}
\end{align}
\end{proposition}

\section{Proofs of Theorem \ref{T1} and Theorem \ref{T2}}
In order to prove our theorems, we first prove the following two conclusions.
\begin{proposition}\label{P2}
Let $M$ be a real hypersurface in $\widetilde{M}^2(c)$ with a $^*$-Ricci soliton whose
potential vector field $W\in\Gamma(T_{\chi}),\chi\neq0$. If the principal curvatures are constant along $\xi$ and $A\xi$ then $M$ is Hopf.
\end{proposition}
\proof Suppose that $M$ is not Hopf, then $\mathcal{N}$ is not empty. Write $W=a_1e_1+a_2e_2+a_3e_3$, where $a_1,a_2,a_3$ are the smooth functions on $\mathcal{N}$. Since $\chi\neq0$, $a_3=0$ and $\chi=\lambda_1$ or $\lambda_2$.
Since $a_1,a_2$ are not all zero, without loss of generality, we may assume $a_1\neq0$, then
 \begin{equation*}
  AW=\chi W\Rightarrow
                     \chi=\lambda_1\quad\hbox{and}\quad
                     a_2=0\quad\hbox{since}\quad\lambda_1\neq\lambda_2.
\end{equation*}
Thus the potential vector field can be written as
$$W=a_1\cos\theta\xi+a_1\sin\theta U.$$
Replacing $f_1$ in \eqref{22} and $f_2$ in \eqref{23} by $a_1\cos\theta$ and $a_1\sin\theta$, respectively, we have
\begin{equation}\label{29}
  \xi(a_1\cos\theta)=\lambda,\quad U(a_1\sin\theta)=0
\end{equation}
because $c=\lambda$ followed from \eqref{24}.
Similarly, in view of the first equation of \eqref{25}, we obtain
\begin{equation}\label{30}
  \xi(a_1\sin\theta)+U(a_1\cos\theta)=0.
\end{equation}

 With the help of \eqref{29} and \eqref{30},  we further obtain
 $$a_1(\sin\theta\xi(\theta)-\cos\theta U(\theta))=-\lambda\sin^2\theta.$$
By \eqref{27}, $\lambda\sin^2\theta=0$. If $\sin\theta\neq0$ then $\lambda=0.$
This leads to a contradiction because $\lambda=c\neq0$. If $\sin\theta=0$ then $W=a_1\cos\theta\xi$, i.e. $\xi$ is a principal vector, which is also a contradiction.
Therefore we complete the proof.\qed\bigskip

\begin{proposition}\label{P1}
A real hypersurface in $\mathbb{C}P^2$, admitting a $^*$-Ricci soliton whose
potential vector field $W\in\Gamma(T_0)$, is Hopf.
\end{proposition}
\proof Suppose that $M$ is not Hopf, then $\mathcal{N}$ is not empty. We may write $W=b_1e_1+b_2e_2+b_3e_3$ in the basis $\{e_1,e_2,e_3\}$,
where $b_1,b_2,b_3$ are smooth functions on $\mathcal{N}.$
By Lemma \ref{L4.1}, $A\phi U=0$, so $AW=0$ implies
 $b_1=b_2=0$, i.e. $W=b_3\phi U$ with $b_3\neq0$. Hence \eqref{25} becomes
\begin{equation}\label{31}
k_3=-\gamma,\quad \xi(b_3)=0,\quad U(b_3)=0.
\end{equation}
And the formula \eqref{23} becomes
\begin{equation}\label{32}
  -b_3\gamma=c-\lambda.
\end{equation}

Since $b_1=0$, the relation \eqref{22} becomes
\begin{equation}\label{33}
-b_3\beta=\lambda.
\end{equation}
So by taking the differentiation of the formula \eqref{32} along $\phi U$, we derive from \eqref{24} that
\begin{equation}\label{34}
  b_3\phi U(\beta)=(c-\lambda)\beta.
\end{equation}
On the other hand, it follows from \eqref{32} and \eqref{33} that
\begin{equation}\label{35}
  \frac{\gamma}{\beta}=\frac{c}{\lambda}-1.
\end{equation}

If $c=\lambda$, then the equation \eqref{35} shows $\gamma=0$.   Further, in view of \eqref{24} we find $\phi U(b_3)=\lambda-c=0$, which means that $b_3$ is constant
since $\xi(b_3)=U(b_3)=0$. Now we derive from \eqref{33} that $\beta$ is constant. Hence together \eqref{17} with \eqref{18}, we obtain
$\beta^2=-\frac{c}{4}.$ It is impossible.

Next we assume $c\neq\lambda$. Thus the relation \eqref{35} follows $\gamma\neq0$ and the formula \eqref{15} follows from \eqref{31}
\begin{equation*}
  \phi U(\beta)=\beta^2-\gamma^2+\frac{c}{4}.
\end{equation*}
Substituting this into \eqref{34}, we get from \eqref{32} that
\begin{equation*}
  (\beta^2-\gamma^2+\frac{c}{4})\frac{1}{\gamma}=-\beta\Rightarrow 1-(\frac{\gamma}{\beta})^2+\frac{c}{4\beta^2}=-\frac{\gamma}{\beta},
\end{equation*}
which reduces from \eqref{35} that $\beta$ is constant. Finally we derive a contradiction from \eqref{34}.
Hence we complete the proof of proposition.\qed\bigskip

{\it Proof of Theorem \ref{T1}.} Under the hypothesis of Theorem \ref{T1}, by Proposition \ref{P2}, $M$ is a Hopf hypersurface of $\widetilde{M}^2(c)$, i.e. $A\xi=\alpha\xi$.
Due to \cite[Theorem 2.1]{NR}, $\alpha$ is constant. We consider a point $p\in M$ and a unit vector field $e\in\mathcal{D}_p$ such
that $Ae=\kappa e$ and $A\phi e=\nu\phi e$, where $\kappa,\nu$ are smooth functions on $M$. Then $\{\xi,e,\phi e\}$ is a local orthonormal basis of $M$. By Corollary 2.3 in \cite{NR},
\begin{equation}\label{36}
  \kappa\nu=\frac{\kappa+\nu}{2}\alpha+\frac{c}{4}.
\end{equation}
Moreover, by a straightforward computation, we have the following lemma.
\begin{lemma}\label{L3}
With respect to $\{\xi,e,\phi e\}$ the Levi-Civita connection is given by
\begin{equation*}\begin{aligned}
  &\nabla_e\xi=\kappa\phi e,\quad\quad\nabla_{\phi e}\xi=-\nu e,\quad\nabla_\xi\xi=0, \\
 &\nabla_ee=a_1\phi e,\quad\quad\nabla_{\phi e}e=\nu\xi+a_2\phi e,\quad\nabla_\xi e=a_3\phi e, \\
 &\nabla_e\phi e=-a_1 e-\kappa\xi,\quad\nabla_{\phi e}\phi e=-a_2e,\quad\nabla_\xi\phi e=-a_3e,
\end{aligned}
\end{equation*}
where $a_1= g(\nabla_ee,\phi e),\, a_2=g(\nabla_{\phi e}e,\phi e),\,a_3=g(\nabla_\xi e,\phi e)$
 are smooth functions on $M$.
\end{lemma}
Under the orthonormal basis $\{\xi,e,\phi e\}$ we may assume that there are
smooth functions $g_1,g_2,g_3$ such that the potential vector filed $W$ can be written as
\begin{equation*}
  W=g_1\xi+g_2e+g_3\phi e.
\end{equation*}
Since $AW=\chi W$ with $\chi\neq0$, we get $\alpha g_1=\chi g_1$, $\kappa g_2=\chi g_2$ and $\nu g_3=\chi g_3.$

Next we consider the following cases:\\

$\bullet$ Case I: $g_1,g_2,g_3$ are not equal zero.\bigskip\\
Then $\kappa=\nu=\alpha$, which leads to $c=0$ from \eqref{36}.
This is a contradiction.\bigskip\\

$\bullet$ Case II: Only one of $g_1,g_2,g_3$ is equal zero.\bigskip\\
If $g_1=0$, then $\kappa=\nu$.
  The relation \eqref{36} yields $(\kappa-\frac{\alpha}{2})^2=\frac{\alpha^2+c}{4}$,
  which shows $\kappa=\nu=constant$ and $\alpha\neq\kappa$; If $g_2=0$, then $\alpha=\nu$,
  \eqref{36} implies $\kappa=\frac{c+2\alpha^2}{2\alpha}$ with $\kappa\neq\alpha$; If $g_3=0$, then $\kappa=\alpha$, which implies $\nu=\frac{c+2\alpha^2}{2\alpha},\nu\neq\alpha$ by \eqref{36}.\bigskip\\

$\bullet$ Case III: Two of $g_1,g_2,g_3$ are equal zero.\bigskip\\
 When $g_1=g_2=0.$ The formula \eqref{10} for $X=\xi$ and $Y=e$ implies
  \begin{equation*}
    g(\nabla_\xi W,e)+g(\xi,\nabla_eW)=0.
  \end{equation*}
  In view of Lemma \ref{L3}, a simple calculation leads to $\kappa=-a_3.$ On the other hand, the relation \eqref{8} for $X=e$ and $Y=\xi$
  yields
  $(\nabla_eA)\xi-(\nabla_\xi A)e=-\frac{c}{4}\phi e.$
  By Lemma \ref{L3}, we find
  \begin{equation}\label{37}
  \alpha\kappa-\kappa\nu-\kappa a_3+a_3\nu=-\frac{c}{4}.
\end{equation}
  A similar computation using relation \eqref{8} for $X=\phi e, Y=\xi$ yields
  \begin{equation}\label{38}
  -\alpha\nu+\kappa\nu-\kappa a_3+a_3\nu=\frac{c}{4}.
\end{equation}
  Moreover, inserting $\kappa=-a_3$ into the above equation gives
\begin{equation}\label{39}
  \kappa^2-\alpha\nu=\frac{c}{4}.
\end{equation}

The combination of \eqref{37} and \eqref{38} leads to $(\kappa-\nu)(2\kappa+\alpha)=0$ because $a_3=-\kappa.$
 If $\nu=\kappa$ then $\alpha\neq\kappa$, otherwise, the formula \eqref{39} will lead to $c=0$.
 If $\nu\neq\kappa$ then $\kappa=-\frac{\alpha}{2}$ and $\nu=\frac{\alpha^2-c}{4\alpha}$.

   When $g_1=g_3=0$, we put $X=\xi, Y=\phi e$ in \eqref{10}. By Lemma \ref{L3}, $a_3=-\nu$, so we get
  $(\kappa-\nu)(2\nu+\alpha)=0$ from \eqref{37} and \eqref{38}. If $\kappa=\nu$ then $\alpha\neq\nu$ as before.
  If $\kappa\neq\nu$ then $\nu=-\frac{\alpha}{2}$ and $\kappa=\frac{\alpha^2-c}{4\alpha}.$

  When $g_2=g_3=0$ the relation \eqref{8} for $X=e, Y=\phi e$ leads to $c=0$ by Lemma \ref{L3}, which is a contradiction.

In a word we have proved that there are two or three distinct constant principal curvatures on $M$.
For the case of $\mathbb{C}P^2$, by Theorem \ref{B} and \cite[Theorem 4.1]{W}, $M$ is an open part of a hypersphere, or a tuber around the complex quadric.

For the case of $\mathbb{C}H^2$, if $M$ has three distinct principal curvatures, by the proof of \cite{BD}, we know that the ruled real hypersurfaces can not be Hopf, which is a contradiction with Proposition \ref{P2}. Thus in this case $M$ has only two distinct constant principle curvatures. In view of Theorem \ref{B}, the real hypersurface $M$  is one of Type $A_{11},A_{2},B$ and $N$.

This finishes the proof of Theorem \ref{T1}.\qed\bigskip

{\it Proof of Theorem \ref{T2}.} Under the assumption of Theorem \ref{T2}, by Proposition \ref{P1} we know that $M$ is a Hopf hypersurface of $\mathbb{C}P^2$.
Hence the equation \eqref{36} and Lemma \ref{L3} are valid. We adopt the same notations as the proof of Theorem \ref{T1}.

Since $AW=0$, we have $\alpha g_1=\kappa g_2=\nu g_3=0.$ If $\alpha=0$
then it follows from \eqref{36} that $\kappa\nu=\frac{c}{4}$, which means that $\kappa,\nu$ are non-zero. So we get $g_2=g_3=0$.
From the Case III in the proof of Theorem \ref{T1}, we know it is impossible.

In the following we assume $\alpha\neq0$, then $g_1=0$. If $g_2$ is also equal zero, then $g_3$ must be non-zero, and further we obtain $\nu=0$ and $\kappa=-\frac{\alpha}{2}\neq0$ from the
Case III in the proof of Theorem \ref{T1}.
If $g_2$ is non-zero then $\kappa=0.$ The formula \eqref{36} implies $\alpha\nu=-\frac{c}{2}$, that shows $\nu$ is a non-zero constant.
Further we know $\alpha\neq\nu$ since $c>0$.

Summarizing the above discussion, we have proved that there are three distinct constant principal curvatures in $M$. Therefore we complete the proof of Theorem \ref{T2} by \cite[Theorem 4.1]{W}.\qed

\section{Proof Theorem \ref{T3}}
In this section we suppose that $M$ is a real hypersurface of $\mathbb{C}P^2$ with a $^*$-Ricci soliton whose potential vector field $W$ belongs to the holomorphic distribution $\mathcal{D}.$ First we prove the following result:
\begin{proposition}\label{P3}
Let $M$ be a real hypersurface in $\mathbb{C}P^2$ with a $^*$-Ricci soliton whose
potential vector field $W\in\mathcal{D}$. If the principal curvatures are constant along $\xi$ and $A\xi$ then $M$ is Hopf.
\end{proposition}
\proof If $M$ is not Hopf then $\mathcal{N}$ is not empty. Let $W=c_1e_1+c_2e_2+c_3e_3\in\mathcal{D}$, where $c_i$ are smooth functions on $\mathcal{N}$, then
\begin{equation}\label{40}
c_1\cos\theta+c_2\sin\theta=0.
\end{equation}
The formula \eqref{22} becomes
\begin{equation}\label{41}
  -c_3\beta=\lambda.
\end{equation}
And by Proposition \ref{P3.2},  the equations \eqref{23}-\eqref{25} accordingly become
\begin{align}
  U(c_1)\sin\theta-U(c_2)\cos\theta-c_3k_3+c-\lambda=0,\label{42} \\
  \phi U(c_3)+c-\lambda=0\label{43}
\end{align}
and
\begin{equation}\label{44}
\left\{
  \begin{array}{ll}
  \xi(c_1)\sin\theta-\xi(c_2)\cos\theta-c_3k_3-c_3\gamma =0, \\
(c_1\sin\theta-c_2\cos\theta)k_3+\xi(c_3)=0,\\
(c_1\sin\theta-c_2\cos\theta)k_1+U(c_3)+\phi U(c_1\sin\theta-c_2\cos\theta)=0.
  \end{array}
\right.
\end{equation}
If $c_3=0$, then the equations \eqref{41} and \eqref{43} show $c=\lambda=0$. It is impossible. Thus $c_3\neq0$, which further implies $\lambda\neq0$ from \eqref{41}.  By \eqref{15} and \eqref{43}, differentiating \eqref{41} along $\phi U$  gives
\begin{equation}\label{45}
    k_3\gamma+\beta^2\frac{c}{\lambda}+\frac{c}{4}=0.
\end{equation}
When $\gamma=0$, this shows $\beta$ is constant. So  it follows from the formula \eqref{15} that $\beta^2=-\frac{c}{4}$, which is impossible because $c>0$.
Hence $\gamma\neq0$ and we get from \eqref{45} that
\begin{equation*}
k_3=-\frac{\beta^2\frac{c}{\lambda}+\frac{c}{4}}{\gamma}
\end{equation*}

If $c_1=c_2=0$, as the proof of Proposition \ref{P1}, by using \eqref{41}-\eqref{44}, we arrive at a contradiction.
Thus one of $c_1,c_2$ must be not zero.

Without loss of generality we set $c_1\neq0$. Taking the differentiation of \eqref{41} along $\xi$ and $U$, respectively, we
obtain from \eqref{28} that $\xi(c_3)=U(c_3)=0$ since $\beta\neq0$. In view of the second equation of \eqref{44} and
\eqref{40}, we find $k_3=0,$ that is,
\begin{equation*}
\beta^2\frac{c}{\lambda}+\frac{c}{4}=0,
\end{equation*}
thus $\beta$ is constant. As before from \eqref{15} we have $\beta^2=-\frac{c}{4}$,
which is impossible. This finishes the proof.\qed\bigskip

{\it Proof of Theorem \ref{T3}}. Under the hypothesis of Theorem \ref{T3}, by Proposition \ref{P3} we know that $M$ is a Hopf hypersurface of $\mathbb{C}P^2$. That means that the structure vector field $\xi$
is a principal vector field, i.e. $A\xi=a\xi$ and $a$ is constant as before.

For any point $p\in M$ we consider a unit vector $Z\in\mathcal{D}_p$
such that $AZ=\mu Z$, then
the following relation holds (see \cite[Corollary 2.3]{NR}):
\begin{equation*}
(\mu-\frac{a}{2})A\phi Z=(\frac{\mu a}{2}+\frac{c}{4})\phi Z.
\end{equation*}

If $\mu=\frac{a}{2}$ the above equation implies $\frac{\mu a}{2}+\frac{c}{4}=0$, i.e. $\mu^2+\frac{c}{4}=0,$ that is impossible. Hence
 $\mu\neq\frac{a}{2}$, which means that $\phi Z$ is a principal vector with principal curvature $\nu$ satisfying
\begin{equation}\label{46}
  \mu\nu=\frac{\mu+\nu}{2}a+\frac{c}{4}.
\end{equation}
Now we know that ${\rm Span}\{Z,\phi Z\}=\mathcal{D}_p$ and $\{\xi,Z,\phi Z\}$ is an orthonormal basis of $T_pM$.
By a straightforward computation, we have
\begin{equation*}
  \nabla_Z\phi Z=-g(\nabla_ZZ,\phi Z)Z-\mu\xi,\quad\nabla_{\phi Z}Z=\nu\xi+g(\nabla_{\phi Z}Z,\phi Z)\phi Z.
\end{equation*}
Taking $X=Z$ and $Y=\phi Z$ in \eqref{8} and using the above formulas, we get
\begin{equation*}
  \mu\nu-\nu a=\frac{c}{4}.
\end{equation*}
Next we divide into two cases to discuss.

  \item[{\bf Case 1.}] If $a\neq0$ then it follows $\mu=\nu$ by combining with \eqref{46} and further $\mu,\nu$ are constant.
Furthermore, we find $\mu=\nu\neq a$, otherwise, the above formula will lead to $c=0$.  By Theorem \ref{B} we get that $M$ is of Type $A_1$.
  \item[{\bf Case 2.}] We assume $a=0$, then $\mu\nu=\frac{c}{4}$. In this case $M$ is a $^*$-Einstein hypersurface(see \cite[Remark 1]{IR}).
The $^*$-Ricci soliton equation \eqref{1} shows $W$ is a conformal Killing vector field, i.e. $\mathcal{L}_Wg=2(\lambda-5c)g$.
From \eqref{7}, we calculate the Ricci operator
\begin{equation*}
  QX=\frac{c}{4}\{5X-3\eta(X)\xi\}+hAX-A^2X,\quad\text{for all}\; X\in TM,
\end{equation*}
where $h=\mathrm{trace}(A)$. Hence by a direct computation we can get that the scalar curvature $r=3c+2\mu\nu$.

Notice that on an $n$-dimensional Riemannian manifold a conformal Killing vector field $X$, i.e. $\mathcal{L}_Xg=2\rho g$, satisfies
\begin{equation*}
  \mathcal{L}_Xr=2(n-1)\Delta\rho-2\rho r,
\end{equation*}
where $r$ is the scalar curvature (see Eq.(5.38) in \cite{Y}). Since $\mu\nu=\frac{c}{4}$, the scalar curvature $r=\frac{7c}{2}\neq0$. Using the above formula we find that $W$ is a Killing vector field.

Moreover, since $M$ is $^*$-Einstein, we derive from Theorem \ref{C} that $M$ is one of Type $A_1,A_2$, and $B$.
But according to the list of principal curvatures of Type $A_1,A_2$ and $B$ hypersurfaces(see \cite[Theorem 3.13-3.15]{NR}), we find that in this case only Type $A_1$  is satisfied.

Therefore we complete the proof of Theorem \ref{T3}.\qed

\section{Proof of Theorem \ref{T4} }

Since the tangent bundle $TM$ can be decomposed as $TM=\mathbb{R}\xi\oplus\mathcal{D}$, where
$\mathcal{D}=\{X\in TM:\eta(X)=0\}$. Then $A\xi$ can be written as
\begin{equation}\label{47}
  A\xi=a\xi+V,
\end{equation}
where $V\in\mathcal{D}$ and $a$ is a smooth function on $M$. In this
section we assume that the hypersurface $M$ of $\widetilde{M}^n(c)$ is equipped with
a $^*$-Ricci soliton such that the potential vector field $W=V$.
\begin{lemma}\label{L1}
On $M$ the following equation is valid:
\begin{equation}\label{48}
     (\nabla_\xi A)\xi=Da+2A\phi V,
\end{equation}
where $Da$ denotes the gradient vector field of $a$.
\end{lemma}
\proof By \eqref{6} and \eqref{47}, for any vector field $X$
\begin{equation}\label{49}
\begin{aligned}
  (\nabla_XA)\xi&=\nabla_X(A\xi)-A\nabla_X\xi \\
   & =X(a)\xi+a\nabla_X\xi+\nabla_XV-A\phi AX.
\end{aligned}
\end{equation}
Thus
\begin{align*}
  g((\nabla_XA)\xi,\xi) & =X(a)+g(\nabla_XV,\xi)-g(A\phi AX,\xi) \\
   &=X(a)-g(V,\nabla_X\xi)-g(\phi AX,A\xi)\\
   &=X(a)+2g(AX,\phi V).
\end{align*}
From the well-known relation $g((\nabla_XA)\xi,\xi)=g((\nabla_\xi A)\xi,X)$ (see \cite[Corollary 2.1]{NR}),  we arrive at \eqref{48}.\qed\bigskip

Next it follows from \eqref{49} and \eqref{8} that
\begin{equation}\label{50}
\begin{array}{ll}
  \nabla_XV & =(\nabla_XA)\xi-X(a)\xi-a\nabla_X\xi+A\phi AX \\
   & =(\nabla_\xi A)X-\frac{c}{4}\phi X-X(a)\xi-a\phi AX+A\phi AX.
\end{array}
\end{equation}
Therefore, by Lemma \ref{L1} we have
\begin{equation}\label{51}
\begin{aligned}
  \nabla_\xi V=&(\nabla_\xi A)\xi-\xi(a)\xi-a\phi A\xi+A\phi A\xi\\
  =&-a\phi V+Da-\xi(a)\xi+3A\phi V.
\end{aligned}
\end{equation}

 Since $\eta(V)=0$, differentiating this along any vector $X$, we have
\begin{equation}\label{52}
g(\nabla_XV,\xi)+g(V,\phi AX)=0.
\end{equation}
In particular, by taking $X=\xi$ in \eqref{52}, we find $g(\nabla_\xi V,\xi)=0$
because of $\nabla_\xi\xi=\phi V$. Hence, taking into account $X=Y=\xi$ in \eqref{10},
we conclude that $\lambda=0.$

Take $X=\xi$ and $Y=\xi$ respectively in \eqref{10}, and it follows from \eqref{51} and \eqref{52} that
\begin{align*}
  -a\phi V+Da-\xi(a)\xi+4A\phi V-2\phi A\phi V=0,\\
  -a\phi V+Da-\xi(a)\xi+4A\phi V=0.
\end{align*}
Hence $\phi A\phi V=0$, which implies $A\phi V=0$ because of \eqref{3} and \eqref{47}. Differentiating
$A\phi V=0$ along vector field $\xi$ and using the first equation of \eqref{6}, \eqref{50} and \eqref{51}, we get
\begin{align*}
  0=\nabla_\xi(A\phi V) & =(\nabla_\xi A)\phi V+A(\nabla_\xi\phi)V+A\phi(\nabla_\xi V) \\
   & =\nabla_{\phi V}V-\frac{c}{4}V+(\phi V)(a)\xi+aAV+A\phi(Da).
\end{align*}
Therefore
\begin{equation}\label{53}
  \nabla_{\phi V}V=\frac{c}{4}V-(\phi V)(a)\xi-aAV-A\phi(Da).
\end{equation}
If we put $X=Y=\phi V$ in \eqref{10}, then the Eq.\eqref{53} leads
to $nc|V|^2=0,$ i.e. $V$ is a zero vector field. Since $\lambda=0$, the following proposition is proved:

\begin{proposition}\label{p}
Every real hypersurface in a non-flat complex space form $\widetilde{M}^n(c)$, $n\geq2$, admitting a $^*$-Ricci soliton with
 potential vector field $V$, is a $^*$-Ricci flat Hopf hypersurface.
\end{proposition}

{\it Proof Theorem \ref{T4}}.  Let $M$ be a $^*$-Ricci flat Hopf hypersurface, namely, $A\xi=a\xi$ and $Q^*X=0$ for all $X$, where $a$ is constant.
In view of \eqref{9}, we have $\frac{cn}{2}\phi^2X+(\phi A)^2X=0$ for all $X$, which further implies
\begin{equation}\label{54}
  \frac{cn}{2}\phi X+A\phi AX=0.
\end{equation}
For any point $p\in M$, let $Z\in\mathcal{D}_p$ is a principal vector, namely, there is a certain function $\mu_1$ such that $AZ=\mu_1 Z$, then
\begin{equation}\label{55}
\mu_1 A\phi Z=-\frac{cn}{2}\phi Z,
\end{equation}
 which shows that $\phi Z$ is also a principal curvature vector, i.e. $A\phi Z=\nu\phi Z$ with $\nu=-\frac{cn}{2\mu_1}$.
On the other hand, as before we know the following relation is also valid:
\begin{align}\label{56}
  (\mu_1-\frac{a}{2})A\phi Z=(\frac{\mu_1 a}{2}+\frac{c}{4})\phi Z.
\end{align}

In the following we divide into two cases.\\

$\bullet$ Case I: $a^2+c\neq0$.\bigskip\\
If $\mu_1=\frac{a}{2}$ then $\frac{\mu_1 a}{2}+\frac{c}{4}=0$, which is a contradiction. Hence $\mu_1\neq\frac{a}{2}$ and from \eqref{56} we find that
the principal curvature $\nu$ is also equal
$(\frac{\mu_1 a}{2}+\frac{c}{4})\Big/(\mu_1-\frac{a}{2})$. Hence we obtain that $\mu_1$ satisfies
\begin{equation}\label{57}
2a\mu_1^2+(1+2n)c\mu_1-acn=0,
\end{equation}
from which we can see that $\mu_1$ is constant. Thus $M$ has constant principal curvatures. However, since $M$ is $^*$-Ricci flat, in view of Theorem \ref{A} and Section 3 in \cite{HA}, we find that there are no hypersurfaces in $\mathbb{C}P^n$ satisfying this case.

For the case of $\mathbb{C}H^n$, in terms of Section 3 in \cite{HA}, only Type $A_{11}$ and $A_{12}$ hypersurfaces may be $^*$-Ricci flat.
 But for the Type $A_{12}$,  we further get $2n=\tanh^2(u)$, which is impossible since $0<\tanh(u)<1.$\\

$\bullet$ Case II: $a^2+c=0$.\bigskip\\
In this case the ambient space is $\mathbb{C}H^n$, since $c=-a^2<0$, $a\neq0$. If $\mu_1\neq\frac{a}{2}$, by \eqref{57}, we get
$\mu_1=na$ and $\nu=\frac{a}{2}$. If $\mu_1=\frac{a}{2}$ then $\nu=na$. Hence it is proved that there are three distinct constant principal curvatures for all $p\in M$.

However, since $M$ is a Hopf, in terms of Theorem \ref{A} and the analysis of Section 3 in \cite{HA}, we know that the Type $A_2$ hypersurfaces can not be $^*$-Einstein, and the Type $B$ and Type $N$ hyersurfaces can not be $^*$-Ricci flat.\\

Summarizing this two cases we complete the proof of Theorem \ref{T4}.\qed
%%% Acknowledgment %%%%%%%%
\section*{Acknowledgement}
The author would like to thank the referees for the helpful
suggestions.

%%%% Bibliography  %%%%%%%%%%


\begin{thebibliography}{99}
\bibitem{B}{\sc J. Berndt}, {\it Real hypersurfaces with constant principal curvatures in complex hyperbolic space,} J. Reine Angew. Math. {\bf 395}(1989), 132-141.
\bibitem{BD}{\sc J. Berndt, J. C. D\'Iaz-Ramos}, {\it Real hypersurfaces with constant principal curvatures in the complex hyperbolic plane,} Pro. A.M.S. {\bf135}(10), 2007, 3349-3357.
\bibitem{BD2} {\sc J. Berndt, J. C. D\'Iaz-Ramos}, {\it Real hypersurfaces with constant principal curvatures in the complex hyperbolic spaces,} J. London Math. Soc. {\bf26}(2), 2006, 778-798.
\bibitem{CR}{\sc T. E. Cecil, P. J. Ryan}, {\it Focal set and real hypersurfaces in complex projective spaces,} Trans. Amer. Math. Soc. {\bf 269}(1982), 481-499.
\bibitem{CK}{\sc J. T. Cho, M. Kimura}, {\it Ricci solitons and real hypersurfaces in a complex space form,}
             Tohoku Math. J. {\bf 61}(2009), 205-212.
\bibitem{CK2}{\sc J. T. Cho, M. Kimura}, {\it Ricci solitons of compact real hypersurfaces in K\"ahler manifolds,}
             Math. Nachr. {\bf 284}(2011), 1385-1393.


\bibitem{H}{\sc R. Hamilton}, {\it The Ricci flow on surfaces, mathematics and general relativity} (Santa Cruz, CA, 1986), Contemp.
         Math. {\bf 71} Amer. Math. Soc., Providence, RI. pp.237-262(1988).

\bibitem{HA}{\sc T. Hamada}, {\it Real hypersurfaces of complex space forms in terms of Ricci $^{*}$-tensor}, Tokyo J. Math. {\bf 25}(2002), 473-483.
\bibitem{IR}{\sc T. Ivey, P. J. Ryan}, {\it The $^*$-Ricci tensor for hypersurface in $\mathbb{C}P^n$ and $\mathbb{C}H^n$,} Tokyo J. Math. {\bf 34}(2011), 445-471.
\bibitem{IR2}{\sc T. Ivey, P. J. Ryan}, {\it Hypersurfaces in $\mathbb{C}P^2$ and $\mathbb{C}H^2$ with two distinct principal curvatures}, Glasgow Math. J.
   {\bf58}(2016), 137-152.
\bibitem{KP}{\sc G. Kaimakamis, K. Panagiotidou}, {\it $^*$-Ricci solitons of real hypersurfaces in non-flat complex space forms},
            J. Geom. Phys. {\bf 86}(2014), 408-413.
\bibitem{KI}{\sc M. Kimura}, {\it Real hypersurfaces and complex submanifolds in a complex projective space ,} Trans. Amer. Math. Soc. {\bf 296}(1986), 137-149.

\bibitem{K2}{\sc M. Kon}, {\it On a Hopf hypersurface of a complex space form,}  Diff. Geom. Appl. {\bf 28}(2010), 295-300.



\bibitem{M}{\sc S. Montoel}, {\it Real hypersurfaces of a complex hyperbolic space},  J. Math. Soc. Japan {\bf 35}(1985), 515-535.

\bibitem{NR}{\sc R. Niebergall, P. J. Ryan}, {\it Real hypersurfaces in complex space forms}, Tight and taut submanifolds
(eds. T. E. Cecil and S. S. Chern), Math. Sci. Res. Inst. Publ. {\bf 32}(1997), Cambridge Univ. Press, 233-305.

\bibitem{PX}{\sc K. Panagiotidou, Ph. J. Xenos},  \emph{Real hypersurfaces in $\mathbb{C}P^2$ and $\mathbb{C}H^2$ whose structure Jacobi operator is
Lie $\mathbb{D}$-parallel}, Note Math. {\bf 32}(2012), 89-99.

\bibitem{T1}{\sc R. Tagaki}, {\it Real hypersurfaces in a complex projective space with constant principal curvatures,} J. Math. Soc. Japan {\bf 27}(1975), 43-53.
\bibitem{T2}{\sc R. Tagaki}, {\it Real hypersurfaces in a complex projective space with constant principal curvatures II,} J. Math. Soc. Japan {\bf 27}(1975), 507-516.
\bibitem{T3}{\sc R. Tagaki}, {\it On homogenous real hypersurfaces in a complex projective space}, Osaka J. Math. {\bf 10}(1973), 495-506.
\bibitem{W}{\sc Q. M. Wang}, {\it Real hypersurfaces with constant principal curvatures in complex complex projective spaces (I),}
         Sci. Sinica. Ser.A {\bf 26}(1983), 1017-1024.
\bibitem{Y}{\sc K. Yano}, {\it Integral Formulas in Riemannian Geometry}, Marcel Dekker Inc.(1970).
\end{thebibliography}
\end{document}